\documentclass{elsart}

\usepackage{amssymb}

\newcommand\cal {\frak}

\newcommand\eq[1] {(\ref{#1})}

\newcommand{\bfm}[1]{\mbox{\boldmath ${#1}$}}

\newcommand{\beqa}{\begin{eqnarray}}
\newcommand{\eeqa}[1]{\label{#1}\end{eqnarray}}
\newcommand{\bequ}{\begin{equation}}
\newcommand{\eequ}[1]{\label{#1}\end{equation}}
\newcommand{\Grad}{\nabla}

\newcommand{\ov}[1]{\overline{#1}}

\newcommand{\Ga}{\alpha}
\newcommand{\Gb}{\beta}
\newcommand{\Gd}{\delta}

\newcommand{\Gve}{\varepsilon}

\newcommand{\Gg}{\gamma}

\newcommand{\Gl}{\lambda}

\newcommand{\Go}{\omega}

\newcommand{\Gz}{\zeta}
\newcommand{\GD}{\Delta}

\newcommand{\GO}{\Omega}

\newcommand{\BGn}{\bfm\eta}

\newcommand{\BGv}{\bfm\nu}

\newcommand{\BGx}{\bfm\xi}

\newcommand{\CN}{{\cal N}}

\newcommand{\BCD}{{\bfm{\cal D}}}

\def\Ba{{\bf a}}

\def\Bx{{\bf x}}
\def\By{{\bf y}}
\def\Bz{{\bf z}}

\def\BD{{\bf D}}

\def\BO{{\bf O}}

\def\BR{{\bf R}}

\def\half{{\scriptstyle{1\over 2}}}

\newcommand{\beq}{\begin{equation}}
\newcommand{\eeq}{\end{equation}}
\newcommand{\overliner}{\begin{eqnarray}}
\newcommand{\earr}{\end{eqnarray}}
\newcommand{\beqn}{\begin{equation*}}
\newcommand{\eeqn}{\end{equation*}}
\newcommand{\overlinern}{\begin{eqnarray*}}
\newcommand{\earrn}{\end{eqnarray*}}
\newcommand{\prt}{\partial}

\newcommand{\fr}{\frac}

\renewcommand{\theequation}{\arabic{equation}}

\begin{document}


\begin{frontmatter}
\title{Uniform asymptotic formulae for Green's kernels in regularly and singularly perturbed domains}
\author{
Vladimir Maz'ya$^1$, Alexander~B.~Movchan$^2$ }
\address{\rm $^1$ Department of Mathematical Sciences, University of Liverpool,  Liverpool L69 3BX,
U.K., ~ and \\ Department of Mathematics,  Ohio State University,
231 W 18th Avenue, Columbus,  OH 43210,
 U.S.A.,  and \\
   Department of Mathematics,  Link\"oping University,  SE-581 83 Link\"oping, Sweden \\
   e-mail: vlmaz@mai.liu.se \\
~\\
  $^2$ Department of Mathematical Sciences, University of Liverpool,  Liverpool L69 3BX,
  U.K. \\
e-mail:
abm@liv.ac.uk\\
}
\maketitle
\thispagestyle{empty}

\begin{abstract}
{Asymptotic formulae for Green's kernels $G_\Gve(\Bx, \By)$ of
various boundary value problems for the Laplace operator are
obtained in regularly perturbed domains and certain domains with
small singular perturbations of the boundary, as $\Gve \to 0$. The
main new feature of these asymptotic formulae is their uniformity
with respect to the independent variables $\Bx$ and $\By$.}
\end{abstract}
\end{frontmatter}
%

\par\medskip\centerline{\rule{2cm}{0.2mm}}\medskip
\setcounter{section}{0}

{\bf Introduction.} In 1907, Hadamard was awarded the Prix
Vaillant by the Acad\'emie des Sciences de Paris for his work
\cite{Had}.
Among much else, this memoir includes investigation of Green's
kernels of some boundary value problems, and, in particular, their
dependence on a small regular perturbation of the domain.
Hadamard's asymptotic formulae played a significant role at the
dawn of functional analysis (see \cite{Levy}). The list of
subsequent applications includes extremal problems in the complex
function theory \cite{Julia}, \cite{Barn-Orn-Pearce}, a biharmonic
maximum principle for hyperbolic surfaces
\cite{Hedenmalm-Jakobsson-Shimorin}, shape sensitivity and
optimisation analysis \cite{Fremiot-Sokolowski}, free boundary
problems \cite{Palmerio-Dervieux}, Brownian motion on
hypersurfaces \cite{Kinateder-McDonald}, theory of reproducing
kernels \cite{Englis-Lukassen}, \cite{Komatsu1},  \cite{Komatsu2}.
Analogues of Hadamard's
formulae were obtained for general elliptic boundary value
problems \cite{Fujuwara-Ozawa} and for the heat equation
\cite{Ozawa}.

The issue of
asymptotic formulae for Green's kernels
in singularly perturbed domains attracted less attention. In this
respect, we
mention the papers \cite{Ozawa1}, \cite{Ozawa2},
where certain estimates for Green's functions in domains with
small holes were obtained.

An important drawback of the estimates just mentioned, which is
also
inherent in the classical Hadamard formulae, is the non-uniformity
with respect to the independent variables. The present article
settles the question of uniformity left open in earlier work, both
for regularly and singularly perturbed domains. First, we sharpen
one of Hadamard's classical formulae in the case of a regular
perturbation. Next, we consider several types of singularly
perturbed domains including domains with a finite number of holes,
thin rods and a truncated cone. For these geometries and different
types of boundary conditions, we derive uniform asymptotic
representations of Green's kernels for the Laplacian.
The asymptotic formulae presented here for
singularly perturbed problems involve Green's functions and
other auxiliary
solutions of certain model boundary value problems independent of
$\Gve$.

The techniques used in the proofs of theorems of the main text of
the paper are based on the method of compound asymptotic
expansions discussed in detail in the two-volume monograph
\cite{MNP}. The uniform asymptotic formulae for Green's kernels
can be efficiently used in numerical algorithms for solving
boundary value problems in domains with singularly perturbed
boundaries.

{\bf 1. Uniform Hadamard's type formula.} Let $\GO$ be a planar
domain with compact closure $\ov{\GO}$ and  smooth boundary $\prt
\GO$. Also, let another domain $\GO(\Gve)$, depending on a small
positive parameter $\Gve$, lie inside $\GO$.
By $\Gd_z$ we denote a smooth positive function defined on $\prt
\GO$, and assume that $\Gve \Gd_z$ is the distance between a point
$\Bz \in \prt \GO$
and $\prt \GO_\Gve$. One of the results in Hadamard's paper
\cite{Had} is the following formula, which relates Green's
functions $G$ and $G_\Gve$ for the Dirichlet boundary value
problem for the Laplacian in  $\GO$ and $\GO_\Gve$: \beq
G_\Gve(\Bx, \By) - G(\Bx, \By) + \Gve \int_{\prt \GO} \fr{\prt
G}{\prt \nu_z}(\Bx, \Bz) \fr{\prt G}{\prt \nu_z}(\Bz, \By) \Gd_z d
s_z = O(\Gve^2). \eequ{eq1} This asymptotic relation holds, in
particular, when either $\Bx$ or $\By$ is placed at a positive
distance from $\prt \GO$, independent of $\Gve$. However, one can
see from the simplest example of two concentric disks $\GO$ and
$\GO_\Gve$ that \eq{eq1} may fail when $\Bx$ and $\By$ approach
the same point at $\prt \GO_\Gve$. In the next theorem, we improve
Hadamard's formula by obtaining a uniform asymptotic
representation for $G_\Gve-G$.

{\bf Theorem 1.} {\it Let $\Bx$ and $\By$ be points of
$\ov{\GO}_\Gve$ situated in a sufficiently thin neighbourhood of
$\prt \GO$.
The right-hand side in {\rm  \eq{eq1}} can be replaced by
\beq \fr{1}{4 \pi} \log \fr{|\Bz(\Bx) - \Bz(\By)|^2 + (\rho_x -
\Gve \Gd_{z(x)} + \rho_y -  \Gve \Gd_{z(y)})^2 }{|z(\Bx) -
z(\By)|^2 + (\rho_x + \rho_y)^2} \eequ{eq2}
$$
+ \fr{1}{2 \pi} \fr{\Gve (\Gd_{z(x)} + \Gd_{z(y)} ) ( \rho_x +
\rho_y) }{|z(\Bx) - z(\By)|^2 + (\rho_x + \rho_y)^2} +
O\Big(\fr{\Gve^2 (\rho_x + \rho_y)}{|z(\Bx) - z(\By)|^2 + (\rho_x
+ \rho_y)^2}\Big),
$$
where $\Bz(\Bx)$ is the point of $\prt \GO$ nearest to $\Bx$,
$\rho_x = | \Bx - \Bz(\Bx)|,$ and the notations $\Bz(\By), \rho_y$
have the same meaning
for the point $\By$.}

We note that the principal term in \eq{eq2} plays the role of a
boundary layer, and when $|\Bz(\Bx) - \Bz(\By)| < \mbox{Const} ~
\Gve $ and $\rho_x + \rho_y < \mbox{Const} ~ \Gve $, it has the
order $O(1)$, whereas the remainder term has the order $O(\Gve)$.

\vspace{.1in}

{\bf 2. Dirichlet problem for a domain with a small inclusion.}
Let $\GO$ and $\Go$ be bounded  domains in $\BR^n$.
We assume that $\GO$ and $\Go$ contain the origin $O$ and
introduce the domain $\Go_\Gve=\{ \Bx: \Gve^{-1} \Bx \in \Go \}.$
Without loss of generality, it is assumed that the minimum
distance between the origin and the points of $\prt \GO$ as well
as the maximum distance between the origin and the points of $\prt
\Go$ are equal to $1$. Let $G_\Gve$ be Green's function of the
Dirichlet problem for the Laplace operator in $\GO_\Gve = \GO
\setminus \ov{\Go}_\Gve.$ We use the notation $|S^{n-1}|$ for the
$(n-1)$-dimensional measure of the unit sphere.

{\bf Theorem 2.} {\it Let $n > 2.$ By $G$ and ${\cal G}$ we denote
Green's functions of the Dirichelt problems in $\GO$ and $\BR^n
\setminus \ov{\Go}$, respectively.

Let $H$ be the regular part of
$G$, that is $H(\Bx, \By) = ( (n-2) |S^{n-1}| )^{-1} |\Bx -
\By|^{2-n} - G(\Bx, \By), $ and let $P$ stand for the harmonic
capacitary potential of $\ov{\Go}$. Then
$$
G_\Gve(\Bx, \By) = G(\Bx, \By) + \Gve^{2-n} {\cal G}(\Gve^{-1}
\Bx, \Gve^{-1} \By) - ( (n-2) |S^{n-1}|  )^{-1}|\Bx - \By|^{2-n}
$$
$$
+ H(0, \By) P( \Gve^{-1} \Bx) + H(\Bx, 0) P(\Gve^{-1} \By)  - H(0,0) P(\Gve^{-1} \Bx) P(\Gve^{-1} \By)
$$
\beq - \Gve^{n-2} ~\mbox{\rm cap}~ \ov{\Go} ~ H(\Bx, 0) H(0, \By)
+
O\Big(\fr{\Gve^{n-1}}{(\min\{|\Bx|, |\By| \})^{n-2}}\Big)
\eequ{eq3} uniformly with respect to $\Bx$ and $\By$ in
$\GO_\Gve$. (Note that the remainder term in \eq{eq3} is $O(\Gve)$
on $\prt
\Go_\Gve$
and $O(\Gve^{n-1})$ on $\prt \GO$. )}

The next theorem contains a result of the same nature for $n=2$.
As before, $G$ and ${\cal G}$ are Green's functions for $\GO$ and
$\BR^2 \setminus \ov{\Go}$, respectively, whereas $H$ is the
regular part of $G$.

{\bf Theorem 3.} {\it Let
$$
\zeta(\BGn) = \lim_{|\BGx| \to\infty} {\cal G}(\BGx, \BGn)
~~\mbox{and}~~ \zeta_\infty = \lim_{|\BGn| \to \infty} \{
\zeta(\BGn) - (2 \pi)^{-1} \log |\BGn|  \}.
$$
Then the asymptotic representation, uniform with respect to $\Bx,
\By \in \GO_\Gve$,  holds
$$
G_\Gve(\Bx, \By) = G(\Bx, \By)  + {\cal G}(\Gve^{-1} \Bx,
\Gve^{-1}\By) + (2 \pi )^{-1} \log(\Gve^{-1} |\Bx - \By| )
$$
$$
+\fr{\Big( (2 \pi)^{-1}\log \Gve + \zeta(\fr{\Bx}{\Gve})
-\zeta_\infty+ H(\Bx, 0) \Big) \Big( (2 \pi)^{-1}\log \Gve +
\zeta(\fr{\By}{\Gve}) -\zeta_\infty + H(0, \By) \Big)}{(2
\pi)^{-1} \log \Gve + H(0,0) - \zeta_\infty}
$$
\beq - \zeta(\Gve^{-1} \Bx) - \zeta(\Gve^{-1} \By) + \zeta_\infty
+O(\Gve). \eequ{eq4}

}

The next assertion is a direct consequence of
\eq{eq4}. It shows that asymptotic representation of $G_\Gve(\Bx,
\By)$ is simplified if $\Bx$ and $\By$ are subject to additional
constraints. We use the same notations as in Theorem 3.


{\bf Corollary.} The following assertions hold:

{\it (a) Let $\Bx$ and $\By$ be points of $\GO_\Gve \subset \BR^n$
such that $\min\{|\Bx|, |\By|\} > 2 \Gve$.
Then for $n=2,$
$$
G_\Gve(\Bx, \By) - G(\Bx, \By) - \fr{G(\Bx,0) G(0, \By)}{({2
\pi})^{-1} \log \Gve + H(0,0) - \Gz_\infty} =
O\big( \fr{\Gve}{\min \{ |\Bx|,|\By|  \}}\big),
$$
and for $n
> 2,$
$$
G_\Gve(\Bx, \By) - G(\Bx, \By) + \Gve^{n-2} \mbox{\rm cap}~
\ov{\Go} ~
G(\Bx, 0) G(0, \By)
= O\big(\fr{\Gve^{n -1}}{ ( |\Bx| |\By|)^{n-1} \min \{ |\Bx|,|\By|
\} }\big).
$$

(b)
If $\max\{|\Bx|, |\By|\} < 1/2$,
then for $n=2,$
$$
G_\Gve(\Bx, \By) - {\cal G}(\fr{\Bx}{\Gve}, \fr{\By}{\Gve}) -
\fr{\Gz(\Gve^{-1}\Bx) \Gz(\Gve^{-1}\By)}{({2 \pi})^{-1} \log \Gve
+ H(0,0) - \Gz_\infty}= O(\max \{|\Bx|, |\By|  \}),
$$
and for $n > 2,$
$$
G_\Gve(\Bx, \By) - \Gve^{2-n}{\cal G}(\fr{\Bx}{\Gve},
\fr{\By}{\Gve}) + H(0,0) (P(\Gve^{-1} \Bx) -1)(P(\Gve^{-1} \By)
-1) = O(\max \{|\Bx|, |\By|  \}).
$$
}

One can say that these formulae are closer, in their spirit, to
Hadamard's original formula \eq{eq1} than Theorems 2 and 3
themselves. Analogous simplified formulae can be deduced easily
from all other theorems of the present paper, but they are not
included here.

\vspace{.1in}

{\bf 3.  Dirichlet-Neumann problems in a planar domain with a
small hole.} The set $\GO_\Gve$ is assumed to be the same as in
Theorem 3, with the additional constraint that $\prt \Go$ is
smooth.
First, let $G_\Gve$ denote the kernel of the inverse operator of
the mixed boundary value problem in $\GO_\Gve$ for the operator
$-\GD$, with the Dirichlet data on $\prt \GO$ and the Neumann data
on $\prt \Go_\Gve$. The notations $G, {\cal G}$ and $H$ have the
same meaning as in Theorem 3, and $\CN$ is the Neumann function in
$\BR^2 \setminus \ov{\Go}$. Let $\BCD$ be a vector function,
harmonic in $\BR^2 \setminus \ov{\Go}$, vanishing at infinity and
such that
$\prt \BCD /\prt \nu = \BGv$ on $\prt \Go$. This vector function
appears in the asymptotic representation
$$
\CN(\BGx, \BGn) \sim (2 \pi)^{-1} \log |\BGx|^{-1} +
(\BCD(\BGn)-\BGn) \cdot \Grad ((2 \pi)^{-1} \log |\BGx|^{-1}),
~~\mbox{as}~~ |\BGx| \to \infty.
$$

{\bf Theorem 4.} {\it The kernel $G_\Gve(\Bx, \By)$ of the inverse
operator of the mixed boundary value problem in $\GO_\Gve$ for the
operator $-\GD$, with the Dirichlet data on $\prt \GO$ and the
Neumann data on $\prt \Go_\Gve$, has the following asymptotic
representation, which is uniform with respect to $\Bx, \By \in
\GO_\Gve,$
$$
G_\Gve(\Bx, \By) = G(\Bx, \By) + \CN(\Gve^{-1} \Bx, \Gve^{-1} \By)
+ (2 \pi)^{-1} \log(\Gve^{-1}|\Bx-\By|)
$$
$$
+ \Gve \BCD(\Gve^{-1} \Bx) \cdot \Grad_x H(0, \By) + \Gve \BCD(\Gve^{-1} \By) \cdot \Grad_y H(\Bx, 0)
$$
\beq
-\Gve^2 \BCD(\Gve^{-1} \By) \cdot ((\Grad_x \otimes \Grad_y) H(0,0)) \BCD(\Gve^{-1} \Bx) + O(\Gve^2).
\eequ{eq5}
}
\vspace{.1in}

Next, consider the
mixed boundary value problem in
$\GO_\Gve$ for the Laplace operator, with the Neumann data on the
smooth boundary $\prt \GO$ and the Dirichlet data on $\prt
\Go_\Gve,$
where $\Go_\Gve = \{\Bx: \Gve^{-1} \Bx \in \Go \subset \BR^2 \}$
with $\Go$ being an arbitray bounded domain.
Let $N(\Bx, \By)$ be
the Neumann function in $\GO$, i.e.
$$
\GD N(\Bx, \By) + \Gd(\Bx-\By) = 0, ~~ \Bx, \By \in \GO, ~~
$$
$$
\fr{\prt}{\prt \nu_x} \Big( N(\Bx, \By) +
({2 \pi})^{-1} \log |\Bx|
 \Big) =0, ~~ \Bx \in \prt \GO, \By \in \GO,
$$
and $$ \int_{\prt \GO} N(\Bx, \By) \fr{\prt}{\prt \nu_x} \log
|\Bx| d s_x =0,$$ with the last condition implying the symmetry of
$N(\Bx, \By)$.
The regular part of the
Neumann function is defined by
$$
R(\Bx,\By) =
({2 \pi})^{-1} \log |\Bx - \By|^{-1} - N(\Bx, \By).
$$
Note that $$R(0,\By) =  -
{(2 \pi)^{-2} } \int_{\prt \GO} \log |\Bx| \fr{\prt}{\prt \nu}
\log |\Bx| d s_x.$$ Let $\BD$ be a vector function, harmonic in
$\BR^2 \setminus \ov{\Go}$, bounded at infinity and subject to the
Dirichlet condition $\BD(\BGx) = \BGx, ~~ \BGx \in \prt {\Go}.$
This vector function appears in the asymptotic representation of
Green's function in $\BR^2 \setminus \ov{\Go}$
$$
{\cal G}(\BGx, \BGn) \sim {\cal G}(\infty, \BGn) + (\BD(\BGn) -
\BGn) \cdot \Grad ((2 \pi)^{-1} \log |\BGx|^{-1}), ~~ \mbox{as} ~~
|\BGx| \to \infty.
$$

{\bf Theorem 5.} {\it The function $G_\Gve(\Bx, \By)$ has the
asymptotic behaviour, uniform with respect to $\Bx, \By \in
\GO_\Gve$:
$$
G_\Gve(\Bx, \By) = {\cal G}(\Gve^{-1} \Bx, \Gve^{-1} \By) + N(\Bx,
\By) - (2 \pi)^{-1} \log |\Bx - \By|^{-1} + R(0,0)
$$
$$
+ \Gve \BD(\Gve^{-1} \By) \cdot \Grad_y R(\Bx, 0) + \Gve
\BD(\Gve^{-1} \Bx) \cdot \Grad_x R(0, \By)
$$
\beq - \Gve^2 \BD(\Gve^{-1} \By) \cdot \Big( ( (\Grad_x \otimes
\Grad_y ) R)(0,0)  \Big)  \BD (\Gve^{-1} \Bx) + O(\Gve^2).
\eequ{eq6} }

\vspace{.1in}

{\bf 4. The Dirichlet-Neumann problem in a thin rod.} Let $C$ be
the infinite cylinder $\{(\Bx', x_n): ~ \Bx' \in \Go, ~ x_n \in
\BR \},$ where $\Go$ is a bounded
domain in $\BR^{n-1}$ with smooth boundary. Also let $C^\pm$
denote Lipschitz subdomains of $C$ separated from $\pm \infty$ by
surfaces $\Gg^\pm$, respectively. We introduce a positive number
$a$, the vector $\Ba=(\BO', a),$  where $\BO'$ is the origin of
$\BR^{n-1},$ and the small parameter $\Gve > 0,$.

By $G_\Gve$ we denote the fundamental solution for $-\GD$ in the
thin domain $C_\Gve = \{ \Bx: \Gve^{-1} (\Bx-\Ba) \in C^+, ~
\Gve^{-1} (\Bx+\Ba) \in C^-\}$ subject to zero Neumann condition
on the cylindrical part of $C$ and zero Dirichlet condition on the
remaining part of $\prt C_\Gve$. Similarly, $G^+$ and $G^-$ stand
for the fundamental solutions for $-\GD$ in the domains $C^\pm$,
and satisfy zero Dirichlet condition on $\Gg^\pm$, zero Neumann
condition on $\prt C^\pm \setminus \Gg^\pm$, and are bounded as
$x_n \to \mp \infty$. Let $\Gz^\pm$ be harmonic functions in
$C^\pm$ subject to  the same boundary conditions as $G^\pm$ and
asymptotically equivalent to $\mp |\Go|^{-1} x_n + \Gz^\pm_\infty$
as $x_n \to \mp \infty$, where $\Gz_\infty^\pm$ are constants. We
note that $\Gz^\pm(\By) = \lim_{\Bx \to \infty} G^\pm (\Bx, \By)$.
Finally, by $G_\infty(\Bx, \By)$ we denote Green's function of the
Neumann problem in $C$ such that
$$
G_\infty(\Bx, \By) = - ( 2 |\Go| )^{-1} | x_n - y_n | + O(\mbox{exp}
(-\Ga |x_n - y_n| )) ~~\mbox{as} ~~ |x_n| \to \infty,
$$
where $\Ga$ is a positive constant, and $|\Go|$ is the
$(n-1)$-dimensional measure of $\Go$.

{\bf Theorem 6.} {\it The following asymptotic formula for
$G_\Gve(\Bx, \By)$, uniform with respect to $\Bx, \By \in
\GO_\Gve,$ holds  } {\small
$$
G_\Gve(\Bx, \By) = \Gve^{2-n} \Big\{
G^+(\Gve^{-1}(\Bx-\Ba),\Gve^{-1}(\By-\Ba)) +
G^-(\Gve^{-1}(\Bx+\Ba),\Gve^{-1}(\By+\Ba)) -
G^\infty(\Gve^{-1}\Bx, \Gve^{-1}\By)
$$
$$
-\Gve \fr{(\fr{x_n}{\Gve |\Go|}  - \half (\Gz_\infty^- -
\Gz_\infty^+) +
\Gz^+(\fr{\Bx-\Ba}{\Gve})-\Gz^-(\fr{\Bx+\Ba}{\Gve})
)(\fr{y_n}{\Gve |\Go|} - \half (\Gz_\infty^- - \Gz_\infty^+) +
\Gz^+(\fr{\By-\Ba}{\Gve})-\Gz^-(\fr{\By+\Ba}{\Gve}))}{2 |\Go|^{-1}
a + \Gve(\Gz_\infty^+ + \Gz_\infty^-)}
$$
$$ +\fr{1}{4} \Big( (\Gve |\Go|)^{-1}{2 a} + \Gz_\infty^- +
\Gz_\infty^+ - 2 \sum_\pm \big( \Gz^\pm(\Gve^{-1}(\Bx\mp \Ba)) +
\Gz^\pm(\Gve^{-1}(\By\mp \Ba)) \big)\Big)
$$
\beq + O(\exp(-\Gb/\Gve))\Big\},\eequ{eq7} } {\it where $\Gb$ is a
positive constant independent of $\Gve$.}

\vspace{.1in}

{\bf 5. The Dirichlet problem in a truncated cone.} Let $K$ be an
infinite cone $\{ \Bx: |\Bx| >0, |\Bx|^{-1} \Bx \in \Go \},$ where
$\Go$ is a subdomain of the unit sphere $S^{n-1}$ such that
$S^{n-1} \setminus \Go$ has a positive $(n-1)$-dimensional
harmonic capacity. Also let $K_0$ and $K_\infty$ denote subdomains
of $K$ separated from
the vertex of $K$ and from $\infty$ by surfaces $\Gg$ and
$\Gamma$, respectively.

By $G_\Gve$ and $G_{cone}$ we denote Green's function of the
Dirichlet problem for $-\GD$ in the domains $K_\Gve= \{ \Bx \in
K_0: \Gve^{-1} \Bx \in K_\infty \}$ and $K$, respectively.
Similarly, $G_0$ and $G_\infty$ stand for Green's functions of the
Dirichlet problem for $-\GD$ in the domains $K_0$ an $K_\infty$
satisfying zero Dirichlet boundary conditions on $\prt K_0
\setminus \{O\}$ and $\prt K_\infty$, with the asymptotic
representations
$$
Z_0(\Bx) = |\Bx|^{2-n-\Gl} \Psi(|\Bx|^{-1} \Bx) (1+o(1)) ~~
\mbox{as} ~ |\Bx| \to 0,
$$
and
$$
Z_\infty(\Bx) = |\Bx|^{\Gl} \Psi(|\Bx|^{-1} \Bx) (1+o(1)) ~~
\mbox{as} ~ |\Bx| \to \infty,
$$
where $\Gl$ is a positive number such that $\Gl (\Gl + n -2)$ is
the first eigenvalue of the Dirichlet spectral problem  in $\Go$
for the Beltrami operator on $S^{n-1}$, and $\Psi$ is the
corresponding eigenfunction. Similarly, we introduce
$\Gl_2 > 0$ such that $\Gl_2 (\Gl_2 + n -2)$ is the second
eigenvalue of the same spectral problem.

{\bf Theorem 7.} {\it Let $\BGx = \Gve^{-1} \Bx$ and $\BGn =
\Gve^{-1} \By$.
Green's function $G_\Gve(\Bx, \By)$ has the asymptotic behaviour
$$
G_\Gve(\Bx, \By) = G_0(\Bx, \By) + \Gve^{2-n} G_\infty(\BGx, \BGn
) - G_{cone}(\Bx, \By)
$$
$$
+ \fr{\Gve^\Gl}{2 \Gl + n -2} \Big\{ \Big( |\BGx|^\Gl
\Psi\Big(\fr{\BGx}{|\BGx|}\Big) - Z_\infty(\BGx) \big) \big(
|\By|^{2-n-\Gl} \Psi\Big(\fr{\By}{|\By|}\Big) - Z_0(\By) \Big)
$$
\beq + \Big( |\BGn|^\Gl \Psi\Big(\fr{\BGn}{|\BGn|}\Big) -
Z_\infty(\BGn) \Big) \Big( |\Bx|^{2-n-\Gl}
\Psi\Big(\fr{\Bx}{|\Bx|}\Big) - Z_0(\Bx) \Big) \Big\} +
O(\Gve^{\min(2 \Gl, \Gl_2)}). \eequ{eq8} This representation is
uniform with respect to $\Bx$ and $\By$ in the truncated cone
$K_\Gve$.}

\vspace{.1in}

{\bf 6. Dirichlet problem in a domain containing several
inclusions.} It is straightforward to generalise the results of
sections dealing with a domain containing a small inclusion/void
to the case of a body containing a finite number of inclusions.
As an example, we formulate a generalisation of Theorem 2. Let
$\GO \subset \BR^n, n > 2,$ be a bounded
domain,
and let $\BO^{(1)}, \BO^{(2)},\ldots, \BO^{(N)}$ be interior
points in $\GO$. Small sets $\Go_\Gve^{(j)}, j=1, \ldots, N,$ are
defined by $ \Go_\Gve^{(j)} = \{\Bx: \Gve^{-1} (\Bx - \BO^{(j)})
\in \Go^{(j)} \subset \BR^n \}, $ where $\Go^{(j)}, ~
j=1,\ldots,N,$ are
bounded domains in $\BR^n$, and they contain the origin $\BO$.
Similar to Section 2, it is assumed that the minimum distance
between $\BO^{(j)}, j=1, \ldots, N,$ and the points of  $\prt \GO$
is equal to $1$. Also, it is supposed that the distance between
$\BO$ and the points of $\prt \Go^{(j)}, j=1, \ldots, N,$ does not
exceed  $1$. By $G_\Gve$
we denote Green's function for the Laplacian in the domain $
\GO_\Gve = \GO \setminus {\cup_{j} \ov{\Go}_\Gve^{(j)}}. $

{\bf Theorem 8.} {\it Let  $G$ and ${\cal G}^{(j)}$ stand for
Green's functions of the Dirichlet problems in $\GO$ and $\BR^n
\setminus \ov{\Go}^{(j)}$, respectively. Also let $H$ be the
regular part of
$G$, and $P^{(j)}$ denote  the harmonic capacitary potentials of
the sets $\ov{\Go}^{(j)}$. Then
$$
G_\Gve(\Bx, \By) = G(\Bx, \By) + \Gve^{2-n} \sum_{j=1}^N {\cal
G}^{(j)}(\fr{\Bx-\BO^{(j)}}{\Gve}, \fr{\By-\BO^{(j)}}{\Gve}) -
\fr{N}{(n-2) |S^{n-1}| |\Bx - \By|^{n-2}}
$$
\beq + \sum_{j=1}^N \Big\{H(\BO^{(j)}, \By) P^{(j)}(
\fr{\Bx-\BO^{(j)}}{\Gve} ) + H(\Bx, \BO^{(j)})
P^{(j)}(\fr{\By-\BO^{(j)}}{\Gve}) \eequ{eq3}
$$
- H(\BO^{(j)},\BO^{(j)}) P^{(j)}(\fr{\Bx-\BO^{(j)}}{\Gve})
P^{(j)}(\fr{\By-\BO^{(j)}}{\Gve})- \Gve^{n-2} ~\mbox{\rm cap}~
\ov{\Go}^{(j)} ~ H(\Bx, \BO^{(j)}) H(\BO^{(j)}, \By) \Big\}
$$
$$ +\sum_{j=1}^N \sum_{1 \leq j \leq N, i \neq j} G(\BO^{(j)},
\BO^{(i)}) P^{(j)}(\fr{\Bx-\BO^{(j)}}{\Gve})
P^{(i)}(\fr{\By-\BO^{(i)}}{\Gve})
$$
$$
 + O\Big(\sum_{j=1}^N
\fr{\Gve^{n-1}}
{(\min\{|\Bx-\BO^{(j)}|, |\By-\BO^{(j)}| \})^{n-2}}\Big) $$
uniformly with respect to $\Bx$ and
$\By$ in $\GO_\Gve$. (Note that the remainder term in \eq{eq3} is
$O(\Gve)$ on $\prt (
\cup_j \Go^{(j)}_\Gve)$ and $O(\Gve^{n-1})$ on $\prt \GO$. )}



\begin{thebibliography}{99}



\bibitem{Had} J. Hadamard, {\it Sur le probl\`eme d'analyse relatif \`a l'\'equilibre des
plaques \'elastiques encastr\'ees.} M\'emoire couronn\'e en 1907
par l'Acad\'emie des Sciences {\bf 33}, No 4, 515-629.


\bibitem{Levy} P. Levy, {\it Le\c{c}ons d'analyse fonctionnelle,} Gauthier-Villars, Paris, 1922.


\bibitem{Julia} G. Julia, {\it Sur une \'equation aux d\'eriv\'ees
fonctionnelles analogue \`a l'\'equation de M. Hadamard}, C.R.
Acad. Sci. Paris {\bf 172} (1921), 831-833.


\bibitem{Barn-Orn-Pearce} R.W. Barnard, K. Pearce, C. Campbell,
{\it A survey of applications of the Julia variation}. Ann. Univ.
Mariae Curie-Sklodowska Sect. A {\bf 54} (2000), 1-20.


\bibitem{Hedenmalm-Jakobsson-Shimorin} H.Hedenmalm, S. Jacobsson,
S. Shimorin, {\it A biharmonic maximum principle for hyperbolic
surfaces}. J. Reine. Angew. Math. 550 (2002), 25-75.



\bibitem{Kinateder-McDonald} K.J. Kinateder, P. McDonald,
{\it Hypersurfaces in $\BR^d$ and the variance of exit times for
Brownian motion}. Proc. Amer. Math. Soc. {\bf 125} (1997), No 8,
2453-2462.

\bibitem{Fremiot-Sokolowski} G. Fremiot, J. Sokolowski,
{\it Shape sensitivity analysis of problems with singularities.
Shape optimization and optimal design} (Cambridge 1999), 255-276,
Lecture Notes in Pure and Appl. Math., 216, Dekker, New York,
2001.


\bibitem{Palmerio-Dervieux} B. Palmerio, A. Dervieux,
{\it Hadamard's variational formula for a mixed problem and an
application to a problem related to a Signorini-like variational
inequality}. Numer. Funct. Anal. Optim. {\bf 1} (1979), no. 2,
113--144.



\bibitem{Englis-Lukassen} M. Englis, D. Lukkassen, J. Peetre, L.E. Persson, {\it On the formula
of Jacques-Louis Lions for reproducing kernels of harmonic and
other functions}. J. Reine Angew. Math. {\bf 570} (2004), 89-129.

\bibitem{Komatsu1} G. Komatsu, {\it Hadamard's variational formula for the Bergman kernel}. Proc. Japan
Acad. Ser. A. Math. Sci. {\bf 58} (1982), No 8, 345-348.

\bibitem{Komatsu2} G. Komatsu, {\it Hadamard's variational formula for the Szeg\"o kernel}. Kodai Math. J.
{\bf 8} (1985), No 2, 157-162.

\bibitem{Marshakov-Wiegmann-Zabrodin} A. Marshakov, P. Wiegmann, A. Zabrodin, {\it Integrable structure
of the Dirichlet boundary problem in two dimensions}. Comm. Math.
Phys. {\bf 227} (2002), No 1, 131-153.

\bibitem{Fujuwara-Ozawa} D. Fujiwara, S. Ozawa, {\it The Hadamard variational formula for the Green functions
of some normal elliptic boundary value problems}. Proc. Japan
Acad. Ser. A Math. Sci. {\bf 54} (1978), No 8, 215-220.

\bibitem{Ozawa} S. Ozawa, {\it Perturbation of domains and Green kernels of heat equations}. Proc. Japan Acad. Ser. A
Math. Sci. {\bf 54} (1978), 322-325; Part II, {\bf 55} (1979),
172-175; Part III, {\bf 55} (1979), 227-230.

\bibitem{Ozawa1} S. Ozawa, S. Roppongi, {\it Heat kernel and
singular variation of domains}. Osaka J. Math. {\bf 32} (1995), no
4, 941-957.

\bibitem{Ozawa2} S. Ozawa, {\it Spectra of the Laplacian with small Robin conditional
boundary}. Proc. Japan Acad. Ser. A Math. Sci. {\bf 72} (1996),
no. 3, 53-54.

\bibitem{MNP} V. Maz'ya, S. Nazarov, B. Plamenevskij,
{\it Asymptotic theory of elliptic boundary value problems in
singularly perturbed domains}, Vols. 1-2, Birkh\"auser, 2000.

\end{thebibliography}
\end{document}